\documentclass[11pt]{article}
\usepackage[a4paper,margin=1in]{geometry}
\usepackage{amsmath,amssymb,amsthm,mathtools}
\usepackage{enumitem}
\usepackage{microtype}
\usepackage{indentfirst}
\usepackage[numbers,sort&compress]{natbib}
\usepackage[hidelinks]{hyperref}
\hypersetup{
  pdftitle={A-alpha-Wronskian Vertices of Paths and Controllability of Path-Rooted Products},
  pdfauthor={Xiaoqi Liu and Haiying Shan},
  pdfsubject={Wronskian vertices of paths and controllability of rooted products},
  pdfkeywords={A-alpha matrix, Wronskian vertex, path, rooted product, walk matrix, controllable graph, resultant}
}
\newcommand{\doi}[1]{\url{https://doi.org/#1}}
\newtheorem{theorem}{Theorem}[section]
\newtheorem{lemma}[theorem]{Lemma}
\newtheorem{proposition}[theorem]{Proposition}
\newtheorem{corollary}[theorem]{Corollary}

\theoremstyle{definition}

\theoremstyle{remark}
\newtheorem{remark}[theorem]{Remark}

\newcommand{\Spec}{\operatorname{Spec}}
\newcommand{\Wr}{\operatorname{Wr}}
\newcommand{\R}{\mathbb{R}}
\newcommand{\Q}{\mathbb{Q}}
\newcommand{\Walpha}{W_{\alpha}}
\newcommand{\Res}{\operatorname{Res}}
\newcommand{\diag}{\operatorname{diag}}

\title{\(A_{\alpha}\)-Wronskian Vertices of Paths\\
and Controllability of Path-Rooted Products}
\author{Xiaoqi Liu\textsuperscript{a}, \, Haiying Shan\textsuperscript{a}\thanks{Corresponding author. E-mail:
		\texttt{2111152@tongji.edu.cn} (X.~Liu),
		\texttt{shan\_haiying@tongji.edu.cn} (H.~Shan).}
	\\ \small \textsuperscript{a}School of Mathematical Sciences, Key Laboratory of Intelligent
	Computing and Applications\\ \small (Ministry of Education), Tongji
    University, Shanghai, China}
\date{}

\begin{document}

\maketitle
\begin{abstract}
For a graph \(G\) with adjacency and degree matrices \(A(G)\) and \(D(G)\),
set \(A_{\alpha}(G)=\alpha D(G)+(1-\alpha)A(G)\) for \(0\leq\alpha<1\).
Using Chebyshev polynomials and a theorem of
Conway and Jones, we classify the \(A_{\alpha}\)-Wronskian vertices
of every path \(P_n\) with \(n\geq2\), namely the vertices \(i\) for
which the standard basis vector \(e_i\) is cyclic for \(A_{\alpha}(P_n)\).
The central vertex of an odd-order path is never Wronskian; for every other
vertex, we determine an explicit finite set of algebraic exceptional
parameters. These sets contain no rational parameters other than
\(0\) and \(1/2\), for which the adjacency and signless Laplacian
cases are characterized by gcd criteria.
For even-order paths rooted at either central vertex, we explicitly evaluate
the root resultant and the walk-determinant polynomial in the
rooted-product factorization.
The resulting controllability criterion combines
\(A_{\alpha}\)-controllability of the base graph with at most three
explicit spectral exclusions.
For each parameter, one of two fixed complementary base graphs
yields controllable products for every even-order path and both central roots.
For \(1/2\leq\alpha<1\), algebraic conjugates further simplify
the spectral test for general rooted products.
\end{abstract}



\section{Introduction}

Throughout this paper, all graphs are finite, simple, and undirected,
with nonempty vertex sets. We write \(\overline G\) for the complement of \(G\).
For a graph \(G\), let \(A(G)\) and \(D(G)\) denote its adjacency
matrix and diagonal degree matrix, respectively. Nikiforov~\cite{Nikiforov2017}
introduced the matrix
\[
  A_{\alpha}(G)=\alpha D(G)+(1-\alpha)A(G).
\]
We assume \(0\leq\alpha<1\) throughout. Writing
\(Q(G)=D(G)+A(G)\) for the signless Laplacian matrix, we have
\[
  A_0(G)=A(G)
  \qquad\text{and}\qquad
  A_{1/2}(G)=\frac{1}{2}Q(G).
\]
For background on signless Laplacian spectra,
see~\cite{CvetkovicRowlinsonSimic2007}.
Related work on \(A_{\alpha}\) includes tree
spectra~\cite{NikiforovPastenRojoSoto2017}, characteristic
polynomials~\cite{LiuLiu2018}, and spectral characterizations of graphs,
including paths~\cite{LinLiuXue2019}.

For a real symmetric matrix \(M\) of order \(N\), we use the following notation
for its characteristic polynomial:
\[
  \phi(M;x)=\det(xI-M).
\]
Let \(M_{\hat v}\) be the
principal submatrix obtained by deleting row and column \(v\);
for a graph matrix we write \(M(G)_{\hat v}\).
For \(A_{\alpha}\), the original vertex degrees are retained.
For \(u\in\R^N\), define the controllability matrix of the pair \((M,u)\) by
\[
  W(M,u)=[\,u,Mu,\ldots,M^{N-1}u\,].
\]
The vector \(u\) is cyclic for \(M\) if these columns span \(\R^N\),
equivalently, if \(\det W(M,u)\neq0\). For the all-ones vector
\(\boldsymbol1\), the matrix \(W(M)=W(M,\boldsymbol1)\) is the walk matrix.
For a graph \(G\), its \(A_{\alpha}\)-walk matrix is
\(\Walpha(G)=W(A_{\alpha}(G))\); we call \(G\)
\(A_{\alpha}\)-controllable if \(\det\Walpha(G)\neq0\).
We write \(\Spec(M)\) for the set of distinct eigenvalues of \(M\).

As in our earlier work~\cite{ShanLiu2025}, we call \(v\) an
\(M\)-Wronskian vertex if the following polynomial has no real zero:
\[
  \Wr\bigl(\phi(M;x),\phi(M_{\hat v};x)\bigr)
  =
  \phi(M;x)\phi'(M_{\hat v};x)
  -
  \phi'(M;x)\phi(M_{\hat v};x).
\]
By Lemma~\ref{lem:wronskian-resultant-cyclic}, this condition is
equivalent to
\[
  \gcd\bigl(\phi(M;x),\phi(M_{\hat v};x)\bigr)=1.
\]
Polynomial gcds are taken to be monic.
This coprimality condition is equivalent to the standard basis vector \(e_v\) being
cyclic for \(M\), or to every eigenvector of \(M\) having a nonzero
\(v\)-th coordinate. For the eigenvector-eigenvalue identity,
see~\cite{DentonParkeTaoZhang2022}; for background on controllability,
see~\cite{Godsil2012,FarrugiaSciriha2014,LiuSiemons2022}.
A rooted graph \(H^{(v)}\) is a graph \(H\) with a specified vertex \(v\).
The rooted product \(G\circ H^{(v)}\) is obtained by attaching one copy
of \(H\) to each vertex of \(G\), identifying its root with that vertex.
We call \(G\) the base graph and \(H^{(v)}\) the rooted factor.
The root condition concerns \(e_v\) in the rooted factor, whereas
controllability of the base graph concerns \(\boldsymbol1\).

Label the vertices of \(P_n\) consecutively by \(1,2,\ldots,n\).
Yan and Wang characterized the \(A_0\)-Wronskian vertices of
\(P_n\)~\cite[Proposition~1 and Remark~1]{YanWang2026}. Since
\(A_{1/2}(P_n)=Q(P_n)/2\), the case \(\alpha=1/2\) is a signless
Laplacian problem. We consider the full range \(0\leq\alpha<1\).
The following theorem gives an explicit classification.

\begin{theorem}[Wronskian vertices of paths]\label{thm:main}
Let \(n\geq2\), \(1\leq i\leq n\), and \(0\leq\alpha<1\). If
\(2i=n+1\), then vertex \(i\) of \(P_n\) is not
\(A_{\alpha}\)-Wronskian. Otherwise, set
\[
  s=\min\{i-1,n-i\},\qquad d=n-2s-1.
\]
Thus \(s\) is the number of vertices on the shorter side of \(i\),
and \(d\) is the difference between the numbers of vertices on the
two sides. For
\(1\leq k\leq d-1\) with
\(\sin(sk\pi/d)\neq0\), define
\[
  \tau_k=-\frac{\sin((s+1)k\pi/d)}{\sin(sk\pi/d)}.
\]
Then vertex \(i\) is \(A_{\alpha}\)-Wronskian if and only if
\begin{equation}\label{eq:main-real-criterion}
  \alpha\notin\mathcal F_{n,i}:=
  \left\{
    \frac{\tau_k}{1+\tau_k}:
    1\leq k\leq d-1,\
    \sin(sk\pi/d)\neq0,\ \tau_k\geq0
  \right\}.
\end{equation}
\end{theorem}

The exceptional set \(\mathcal F_{n,i}\) is finite and consists of
real algebraic numbers, that is, roots of nonzero polynomials over \(\Q\).
The next corollary identifies the rational exceptions.

\begin{corollary}[Rational classification]
\label{cor:rational-classification}
Let \(n\geq2\), \(1\leq i\leq n\), and
\(\alpha\in\Q\cap[0,1)\). Vertex \(i\) of \(P_n\) is
\(A_{\alpha}\)-Wronskian if and only if
\[
  \begin{cases}
    \gcd(i,n+1)=1, & \alpha=0,\\[1mm]
    \gcd(2i-1,n)=1, & \alpha=\dfrac12,\\[2mm]
    2i\neq n+1,
      & \alpha\in\Q\cap(0,1)\setminus\left\{\dfrac12\right\}.
  \end{cases}
\]
\end{corollary}

We next study the controllability of rooted products with paths.
Characteristic-polynomial formulas for rooted products were
obtained by Schwenk~\cite{Schwenk1974} and Godsil and
McKay~\cite{GodsilMcKay1978}; see also~\cite{RajkumarPavithra2022}.
Recent work gives walk-determinant formulas for path-rooted products
in the adjacency case~\cite{WangYanMao2024,YanWang2026}.
For paths rooted at an end vertex, Yan, Mao, and Wang~\cite{YanMaoWang2024}
gave a signless Laplacian formula and conjectured a parameterized extension.
End-rooted paths were also used in \(A_{\alpha}\)-spectral
constructions~\cite{TianWuCui2024}.

The general rooted-product controllability criterion is due to
Shan and Liu~\cite[Theorem~6.1]{ShanLiu2025}.
The walk-determinant factorization of Wang, Shen, and
Mao~\cite[Theorem~3]{WangShenMao2026} also applies to \(A_{\alpha}\)
by the argument in Section~\ref{sec:rooted-general}.
When both graphs have at least two vertices, the criterion requires
a controllable base graph, a Wronskian
root, and a walk-determinant polynomial (defined in
Section~\ref{sec:rooted-general}) that is nonzero on the base spectrum.
For paths, making this criterion explicit requires determining both
the admissible roots and the zeros of the walk-determinant polynomial.
Theorem~\ref{thm:main} supplies the root condition for every path.

For either central vertex of \(P_{2r}\),
Theorem~\ref{thm:even-center} explicitly evaluates the path-dependent
factors. The walk-determinant polynomial has a factor \(t^r\) and two
further factors with parity-dependent exponents, yielding at most three
spectral exclusions for the base graph
(Corollary~\ref{cor:even-center-controllable}).

Theorem~\ref{thm:two-seed-cover} also gives a uniform construction:
two fixed complementary six-vertex graphs cover all \(\alpha\in[0,1)\).
For each parameter, the same choice of base works for every \(r\geq1\)
and both central roots. Two is the minimum possible number of such bases.

Section~2 classifies the Wronskian vertices of paths using Chebyshev
polynomials of the second kind and a theorem of Conway and Jones.
Section~3 uses resultants to obtain the central-root formulas and a
base-graph construction for controllable rooted products.
Section~4 gives a further base-graph criterion for
\(1/2\leq\alpha<1\) using positive semidefiniteness and algebraic conjugates.

\section{Wronskian vertices of paths}

\subsection{Characteristic polynomials of the tridiagonal blocks}
\label{sec:path-blocks}

Deleting a vertex of a path leaves two tridiagonal blocks.
Fix \(0\leq\alpha<1\) and set \(\beta=1-\alpha\) and
\(\tau=\alpha/\beta\).
For \(m\geq1\), let \(B_m\) be the \(m\times m\) tridiagonal matrix
\[
  B_m=
  \begin{pmatrix}
    \alpha & \beta \\
    \beta & 2\alpha & \ddots \\
          & \ddots & \ddots & \beta\\
          &        & \beta  & 2\alpha
  \end{pmatrix},
\]
and define
\[
  f_m(x)=\det(xI_m-B_m),\qquad f_0(x)=1.
\]
The entry \(\alpha\) comes from an end vertex of the original path.
Reversing the vertex order moves this entry to the other end of the
block and leaves its characteristic polynomial unchanged.

\noindent\begin{minipage}{\linewidth}
\begin{lemma}\label{lem:deleted}
For \(1\leq i\leq n\),
\[
  \phi\bigl(\bigl(A_{\alpha}(P_n)\bigr)_{\hat i};x\bigr)
  =f_{i-1}(x)f_{n-i}(x).
\]
Moreover, vertex \(i\) is \(A_{\alpha}\)-Wronskian if and only if
\[
  \gcd\bigl(f_{i-1}(x),f_{n-i}(x)\bigr)=1.
\]
\end{lemma}
\end{minipage}

\begin{proof}
Set \(J=A_{\alpha}(P_n)\), \(p(x)=\phi(J;x)\), and
\(q(x)=\phi(J_{\hat i};x)\). The two blocks of \(J_{\hat i}\)
have characteristic polynomials \(f_{i-1}\) and \(f_{n-i}\).
Their direct sum therefore gives \(q=f_{i-1}f_{n-i}\).

To prove the Wronskian criterion, move coordinate \(i\) to the
first position and write
\[
  J=\begin{pmatrix}a_i&b^T\\ b&C\end{pmatrix},
  \qquad C=J_{\hat i},\qquad a_i=\alpha d_{P_n}(i),
\]
where \(d_{P_n}(i)\) is the degree of \(i\), and \(b\) has entry
\(\beta=1-\alpha\) at each neighbor of \(i\) and zero elsewhere.
Thus \(a_i=\alpha\) at an end vertex and \(a_i=2\alpha\) at an internal
vertex. If \(x\) is not an eigenvalue of \(C\), the Schur complement
formula gives
\[
\begin{aligned}
  p(x)
  &=\det\begin{pmatrix}x-a_i&-b^T\\-b&xI-C\end{pmatrix}\\
  &=\det(xI-C)\bigl(x-a_i-b^T(xI-C)^{-1}b\bigr).
\end{aligned}
\]
Choose an orthonormal eigenbasis \(z_1,\ldots,z_{n-1}\) of \(C\)
with each vector supported on one block, and write \(Cz_j=\mu_jz_j\).
Since each nonempty block is tridiagonal with nonzero off-diagonal
entries, its spectrum is simple and every eigenvector is nonzero at
the end vertex adjacent to \(i\). Hence \(c_j=b^Tz_j\neq0\).
Dividing by \(q(x)=\det(xI-C)\) and using the spectral
decomposition of \(C\) gives
\[
  \frac{p(x)}{q(x)}
  =x-a_i-\sum_{j=1}^{n-1}\frac{c_j^2}{x-\mu_j}.
\]
Differentiating this identity and multiplying by \(-q(x)^2\) gives
\[
  \Wr(p,q)(x)
  =-q(x)^2-\sum_{j=1}^{n-1}c_j^2
    \left(\frac{q(x)}{x-\mu_j}\right)^2.
\]
Each quotient is a polynomial, so the identity holds for all \(x\).
If the block spectra are disjoint, every zero of \(q\) is simple.
At such a zero one term in the sum is strictly negative, and away
from the zeros the term \(-q(x)^2\) is strictly negative.
If the blocks share an eigenvalue \(\mu\), it is a double zero of
\(q\), and every term vanishes at \(\mu\).
Thus the Wronskian has no real zero if and only if
\(f_{i-1}\) and \(f_{n-i}\) are coprime.
For \(i=1,n\), one block is empty and \(f_0=1\), so the same
argument applies.
\end{proof}

Expanding the tridiagonal determinant gives \(f_1(x)=x-\alpha\)
and, for \(m\geq2\),
\begin{equation}\label{eq:frecurrence}
  f_m(x)=(x-2\alpha)f_{m-1}(x)-\beta^2f_{m-2}(x).
\end{equation}

Let \(U_m\) denote the Chebyshev polynomial of the second kind,
normalized by
\[
  U_{-1}(t)=0,\qquad U_0(t)=1,\qquad
  U_m(t)=2tU_{m-1}(t)-U_{m-2}(t)\quad(m\geq1).
\]

\begin{lemma}\label{lem:chebyshev}
Set \(z=(x-2\alpha)/\beta\). Then, for every \(m\geq0\),
\begin{equation}\label{eq:chebyshevform}
  f_m(x)=\beta^m
  \left[
    U_m\left(\frac z2\right)
    +\tau U_{m-1}\left(\frac z2\right)
  \right].
\end{equation}
\end{lemma}

\begin{proof}
The formula holds for \(m=0\) by the initial conditions. At \(m=1\),
the right-hand side equals
\[
  \beta(z+\tau)
  =(x-2\alpha)+\alpha=x-\alpha.
\]
For \(m\geq2\), both sides satisfy the recurrence
\eqref{eq:frecurrence}, so induction proves the formula.
\end{proof}

\subsection{Common zeros and the real-parameter classification}

Define
\[
  Y_m(z)=
  U_m\left(\frac z2\right)
  +\tau U_{m-1}\left(\frac z2\right).
\]
The sequence \(Y_m\) satisfies
\[
  Y_0(z)=1,
  \qquad Y_1(z)=z+\tau,
  \qquad
  Y_m(z)=zY_{m-1}(z)-Y_{m-2}(z)\quad(m\geq2).
\]

Path reversal allows us to assume \(1\leq i\leq(n+1)/2\).
Put \(d=n-2i+1\).

\begin{proposition}\label{prop:trigcriterion}
Suppose \(i<(n+1)/2\). Then \(f_{i-1}\) and \(f_{n-i}\) have a common
zero if and only if there exists an integer \(k\), \(1\leq k\leq d-1\),
such that
\begin{equation}\label{eq:sineratio}
  \sin\frac{(i-1)k\pi}{d}\neq0
  \qquad\text{and}\qquad
  \tau=-\frac{\sin(ik\pi/d)}{\sin((i-1)k\pi/d)}.
\end{equation}
\end{proposition}

\begin{proof}
\smallskip
\noindent\emph{Necessity.}
Induction on \(d\), using the recurrence for \(Y_m\), gives
\[
  Y_{n-i}(z)
  =
  U_{d-1}\left(\frac z2\right)Y_i(z)
  -
  U_{d-2}\left(\frac z2\right)Y_{i-1}(z).
\]
Suppose that \(Y_{i-1}(z)=Y_{n-i}(z)=0\).
The recurrence shows that consecutive polynomials \(Y_{i-1}\)
and \(Y_i\) are coprime: a common zero would propagate back to
\(Y_0=1\). Since \(Y_i(z)\neq0\), we obtain
\[
  U_{d-1}\left(\frac z2\right)=0.
\]
In particular, \(d\geq2\), and
\[
  z=2\cos\theta,\qquad
  \theta=\frac{k\pi}{d}
\]
for some \(1\leq k\leq d-1\). Using
\(U_m(\cos\theta)=\sin((m+1)\theta)/\sin\theta\), the equation
\(Y_{i-1}(z)=0\) becomes
\[
  \sin(i\theta)+\tau\sin((i-1)\theta)=0.
\]
If \(\sin((i-1)\theta)=0\), then \(\sin(i\theta)=0\) as well,
so subtracting the two integer multiples of \(\pi\) gives
\(\theta\in\pi\mathbb Z\), contrary to \(0<\theta<\pi\).
Thus \(\sin((i-1)\theta)\neq0\), and division gives \eqref{eq:sineratio}.

\smallskip
\noindent\emph{Sufficiency.}
Suppose that \eqref{eq:sineratio} holds, and set
\(\theta=k\pi/d\) and \(z=2\cos\theta\). Then \(Y_{i-1}(z)=0\), while
\(U_{d-1}(z/2)=0\).
Substituting these equalities into the expression for \(Y_{n-i}\)
gives \(Y_{n-i}(z)=0\).
By Lemma~\ref{lem:chebyshev}, the corresponding value of \(x\)
is a common zero of \(f_{i-1}\) and \(f_{n-i}\).
\end{proof}

We can now prove the real-parameter classification.

\begin{proof}[Proof of Theorem~\ref{thm:main}]
If \(2i=n+1\), then the two polynomials in
Lemma~\ref{lem:deleted} are equal and nonconstant. Hence the unique middle
vertex of an odd-order path is never \(A_{\alpha}\)-Wronskian. We may
therefore assume \(2i\neq n+1\).

Reverse the path if necessary so that \(i<(n+1)/2\).
Then \(s=i-1\) and \(d=n-2i+1\), and the ratio in
Proposition~\ref{prop:trigcriterion} equals \(\tau_k\).
Thus \(i\) is not Wronskian exactly when \(\tau=\tau_k\) for
one of the indices in the theorem. Since
\[
  \tau=\frac{\alpha}{1-\alpha}\geq0,
  \qquad
  \alpha=\frac{\tau}{1+\tau},
\]
this is equivalent to \(\alpha\in\mathcal F_{n,i}\).

There are finitely many indices \(k\), so \(\mathcal F_{n,i}\)
is finite. The sines of rational multiples of \(\pi\) lie in
cyclotomic fields, namely fields \(\Q(\zeta)\) generated by roots of
unity (complex numbers \(\zeta\) satisfying \(\zeta^L=1\) for some
integer \(L\geq1\)). Hence \(\tau_k\) and \(\tau_k/(1+\tau_k)\)
are algebraic.
\end{proof}

For a fixed \(\alpha\), path reversal reduces the tests in
\eqref{eq:main-real-criterion} to \(1\leq i<(n+1)/2\).
For \(2\leq n\leq6\), direct evaluation of
\eqref{eq:main-real-criterion} gives
\(\mathcal F_{n,i}\subseteq\{0,1/2\}\) for all \(i\) with \(2i\neq n+1\).
Among vertices with \(2i\neq n+1\), the smallest path admitting an
irrational exceptional parameter is \(P_7\).
Taking \(i=2\), \(d=4\), and \(k=3\) in
\eqref{eq:sineratio} gives
\[
  \tau=-\frac{\sin(3\pi/2)}{\sin(3\pi/4)}=\sqrt2,
  \qquad
  \alpha=\frac{\sqrt2}{1+\sqrt2}=2-\sqrt2.
\]
Thus vertices \(2\) and \(6\) are not
\(A_{2-\sqrt2}\)-Wronskian.

\subsection{Rational parameters}

For rational \(\alpha\notin\{0,1/2\}\), we need a restriction on
the rational ratios of consecutive sines in Proposition~\ref{prop:trigcriterion}.
Their absolute values are ratios of chord lengths in regular polygons.
We use the following two special cases of
Conway and Jones~\cite[Theorem~7]{ConwayJones1976}.

\begin{theorem}[Conway and Jones, special cases]
\label{thm:conway-jones}
If \(\varphi\) is a rational multiple of \(\pi\) and
\(\cos\varphi\in\Q\), then
\[
  \cos\varphi\in\left\{0,\ \pm\frac12,\ \pm1\right\}.
\]
Moreover, let \(\varphi_1\) and \(\varphi_2\) be distinct rational
multiples of \(\pi\) in \((0,\pi/2)\), with both cosines irrational.
If \(a,b\in\Q\setminus\{0\}\) and
\[
  a\cos\varphi_1+b\cos\varphi_2\in\Q,
\]
then \(|a|=|b|\).
\end{theorem}

\begin{corollary}\label{cor:chords}
Let \(\ell_1\) and \(\ell_2\) be two nonzero chord lengths of a
regular polygon. If \(\ell_1/\ell_2\in\Q\), then either
\(\ell_1=\ell_2\), or the longer chord is a diameter and is twice the
shorter chord.
\end{corollary}

\begin{proof}
Let \(R\) be the circumradius and choose rational multiples
\(0<\theta_1,\theta_2\leq\pi/2\) of \(\pi\) such that
\[
  \ell_j=2R\sin\theta_j \qquad (j=1,2).
\]
For \(q=\ell_1/\ell_2\in\Q_{>0}\), squaring gives
\[
  \cos(2\theta_1)-q^2\cos(2\theta_2)=1-q^2.
\]

Since \(q^2\) and \(1-q^2\) are rational and \(q>0\), the two cosines
are either both rational or both irrational. If both are irrational,
cosine symmetries give
\[
  \cos(2\theta_j)=\varepsilon_j\cos\varphi_j,
  \qquad
  \varepsilon_j\in\{\pm1\},
  \qquad
  0<\varphi_j<\frac{\pi}{2},
\]
where each \(\varphi_j\) is a rational multiple of \(\pi\).
If \(\varphi_1=\varphi_2\), irrationality forces
\(\varepsilon_1=q^2\varepsilon_2\), hence \(q^2=1\);
if they are distinct, the second assertion of
Theorem~\ref{thm:conway-jones} gives the same conclusion.
Hence \(\ell_1=\ell_2\).

If both cosines are rational, the first assertion of
Theorem~\ref{thm:conway-jones} gives
\[
  \cos(2\theta_j)
  \in\left\{0,\ \pm\frac12,\ \pm1\right\}.
\]
Since the chords are nonzero,
\[
  \frac{\ell_j^2}{R^2}
  =2\bigl(1-\cos(2\theta_j)\bigr)
  \in\{1,2,3,4\}.
\]
Since \(q\) is rational, the quotient of these two squared lengths
must be a rational square. For distinct values in \(\{1,2,3,4\}\),
this occurs only for \(1\) and \(4\). Thus the unequal chord lengths
are \(R\) and \(2R\).
\end{proof}

\begin{lemma}\label{lem:consecutive}
Let \(m\) be a positive integer and \(\theta\) a rational multiple of
\(\pi\), and suppose
\(\sin(m\theta)\sin((m+1)\theta)\neq0\). If
\[
  \frac{|\sin((m+1)\theta)|}{|\sin(m\theta)|}\in\Q,
\]
then the ratio is \(1\).
\end{lemma}

\begin{proof}
Write \(\theta=a\pi/N\), with \(a,N\in\mathbb Z\) and \(N\geq2\).
The central angle between consecutive vertices of a regular
\(2N\)-gon is \(\pi/N\). For unit circumradius, a chord spanning
\(k\) steps therefore has length \(2|\sin(k\pi/(2N))|\).
For \(k=2ma\), the angle in this sine is
\(k\pi/(2N)=m\theta\), and similarly for \(k=2(m+1)a\).
Thus the chords joining vertex \(0\) to vertices \(2ma\) and
\(2(m+1)a\), with indices taken
modulo \(2N\), have lengths \(2|\sin(m\theta)|\) and
\(2|\sin((m+1)\theta)|\). By
Corollary~\ref{cor:chords}, an unequal rational ratio would imply
\[
  \bigl\{|\sin(m\theta)|,|\sin((m+1)\theta)|\bigr\}
  =\left\{1,\frac12\right\}.
\]
Let \(j\in\{m,m+1\}\) correspond to the value \(1\), and let \(j'\)
be the other index. Modulo \(\pi\),
\[
  j\theta\equiv\frac{\pi}{2},
  \qquad
  j'\theta\equiv\frac{\pi}{6}
  \ \text{or}\ 
  \frac{5\pi}{6}.
\]
Since \(|j-j'|=1\), subtraction gives
\(\theta\equiv\pm\pi/3\pmod{\pi}\). As \(j\) is an integer,
\(j\theta\) is therefore congruent to one of
\(0,\pi/3,2\pi/3\) modulo \(\pi\), none of which is congruent to
\(\pi/2\). This contradicts \(j\theta\equiv\pi/2\pmod{\pi}\).
Hence the ratio is \(1\).
\end{proof}

For the adjacency case, we use the following result of Yan and Wang.

\begin{lemma}[Yan and Wang
  {\cite[Proposition~1 and Remark~1]{YanWang2026}}]
\label{lem:adjacency-path-wronskian}
Let \(n\geq2\) and \(1\leq i\leq n\). Vertex \(i\) of \(P_n\) is
\(A_0\)-Wronskian if and only if
\[
  \gcd(i,n+1)=1.
\]
\end{lemma}

\begin{proof}[Proof of Corollary~\ref{cor:rational-classification}]
By path reversal, we may assume \(i\leq(n+1)/2\). Set \(d=n-2i+1\).

The case \(d=0\) is excluded by Theorem~\ref{thm:main}.
Here \(\gcd(i,n+1)=i>1\) and \(\gcd(2i-1,n)=n>1\), so all
three criteria fail. We may therefore assume \(d>0\).

\smallskip
\noindent\emph{Case 1: \(\alpha=0\).}
The assertion is Lemma~\ref{lem:adjacency-path-wronskian}.

\smallskip
\noindent\emph{Case 2: \(\alpha=1/2\).}
Since \(\tau=1\), Proposition~\ref{prop:trigcriterion} shows that
\(f_{i-1}\) and \(f_{n-i}\) have a common zero if and only if
\[
  \sin(i\theta)+\sin((i-1)\theta)=0,
  \qquad
  \theta=\frac{k\pi}{d},
  \quad 1\leq k\leq d-1.
\]
Any such solution has \(\sin((i-1)\theta)\neq0\), since otherwise
both sines would vanish, forcing \(\theta\in\pi\mathbb Z\).
Since \(0<\theta<\pi\),
\[
  0
  =
  2\sin\left(\frac{(2i-1)\theta}{2}\right)
   \cos\left(\frac{\theta}{2}\right)
\]
is equivalent to
\[
  2d\mid (2i-1)k.
\]
Put \(g=\gcd(2i-1,d)\). Since \(2i-1\) is odd, \(g\) is odd
and \(\gcd(2i-1,2d)=g\). The divisibility condition is therefore
\(2d/g\mid k\). Such a \(k\) lies in \(1\leq k<d\) exactly when
\(g>2\), which for odd \(g\) is equivalent to \(g>1\).
Hence the two polynomials have a common zero if and only if
\(\gcd(2i-1,d)>1\).

Finally, \(d=n-(2i-1)\), so
\(\gcd(2i-1,d)=\gcd(2i-1,n)\).
Thus vertex \(i\) is \(A_{1/2}\)-Wronskian if and only if
\(\gcd(2i-1,n)=1\).

\smallskip
\noindent\emph{Case 3:
\(\alpha\in\Q\cap(0,1)\setminus\{1/2\}\).}
For \(i=1\), the polynomial \(f_{i-1}=f_0=1\) has no zero.
We may therefore assume \(i\geq2\). Here
\[
  \tau=\frac{\alpha}{1-\alpha}\in\Q_{>0},
  \qquad \tau\neq1.
\]
If \(f_{i-1}\) and \(f_{n-i}\) had a common zero,
Proposition~\ref{prop:trigcriterion} would imply
\[
  \theta=\frac{k\pi}{d},
  \qquad
  \tau
  =
  \frac{|\sin(i\theta)|}{|\sin((i-1)\theta)|}.
\]
The denominator is nonzero, and \(\tau>0\) allows us to take
absolute values in the sine ratio.
Lemma~\ref{lem:consecutive} then gives \(\tau=1\), a contradiction.
Consequently,
\[
  \gcd(f_{i-1},f_{n-i})=1
\]
whenever \(d>0\). By Lemma~\ref{lem:deleted}, vertex \(i\) is
\(A_{\alpha}\)-Wronskian whenever \(2i\neq n+1\).

The three cases, together with path reversal, prove the corollary.
\end{proof}

\section{Controllable rooted products}

We apply the path classification to rooted-product controllability,
starting with the general criterion.

\subsection{The general rooted-product criterion}
\label{sec:rooted-general}

Let \(G\) have order \(n\), and let \(H^{(v)}\) have order \(m\).
Let \(E_{vv}\) denote the matrix with a \(1\) in position \((v,v)\) and
zeros elsewhere. Then
\begin{equation}\label{eq:Aalpha-rooted}
  A_{\alpha}(G\circ H^{(v)})
  =
  I_n\otimes A_{\alpha}(H)
  +A_{\alpha}(G)\otimes E_{vv}.
\end{equation}
Here \(\otimes\) denotes the Kronecker product: for \(S=(s_{ij})\),
the product \(S\otimes T\) is the block matrix \((s_{ij}T)\).
The second summand adds the edges between roots and the
corresponding contributions to their degrees.

For \(M\in\{A,Q,A_{\alpha}\}\) and a graph \(X\), write
\(W_M(X)=W(M(X))\), with \(\Walpha(X)\) as shorthand when
\(M=A_{\alpha}\). We call \(X\) \(M\)-controllable if
\(\det W_M(X)\neq0\).
For nonzero polynomials \(f(x)=a\prod_{i=1}^{r}(x-\xi_i)\) and
\(g(x)=b\prod_{j=1}^{s}(x-\zeta_j)\), their resultant is
\[
  \Res(f,g)
  =a^s b^r
    \prod_{i=1}^{r}\prod_{j=1}^{s}(\xi_i-\zeta_j).
\]
The complex zeros in these products are counted with multiplicity.
In particular, \(\Res(f,g)\neq0\) if and only if
\(\gcd(f,g)=1\).

We first relate the root condition to resultants.

\begin{lemma}
\label{lem:wronskian-resultant-cyclic}
Let \(M\) be a real symmetric matrix of order \(N\), and let
\(v\in\{1,\ldots,N\}\). Set
\[
  p(x)=\phi(M;x),\qquad q(x)=\phi(M_{\hat v};x).
\]
Then the following are equivalent:
\begin{enumerate}
  \item \(\Wr(p,q)\) is nonzero for every real \(x\);
  \item \(\gcd(p,q)=1\), equivalently, \(\Res(p,q)\neq0\);
  \item the standard basis vector \(e_v\) is cyclic for \(M\).
\end{enumerate}
\end{lemma}

\begin{proof}
Since \(p\) and \(q\) have only real zeros, a common zero would
also be a zero of their Wronskian. When \(M\) has simple spectrum
\(\lambda_1,\ldots,\lambda_N\), choose a unit \(\lambda_j\)-eigenvector
\(u_j\) for each \(j\).
Cramer's rule and the spectral decomposition give
\begin{equation}\label{eq:principal-resolvent}
  \frac{q(x)}{p(x)}
  =e_v^T(xI-M)^{-1}e_v
  =\sum_{j=1}^{N}\frac{|u_j(v)|^2}{x-\lambda_j}.
\end{equation}
Multiplying by \(x-\lambda_j\) and letting \(x\to\lambda_j\) gives
\(q(\lambda_j)=|u_j(v)|^2p'(\lambda_j)\).
Now suppose that \(p\) and \(q\) are coprime. By Cauchy's interlacing
theorem, a repeated eigenvalue of \(M\) would also belong to the spectrum
of \(M_{\hat v}\). Thus \(M\) has simple spectrum, and the identity for
\(q(\lambda_j)\) gives \(u_j(v)\neq0\) for every \(j\).
For \(x\) outside the spectrum of \(M\), differentiating yields
\[
  \Wr(p,q)
  =p(x)^2\left(\frac{q}{p}\right)'(x)
  =-p(x)^2\sum_{j=1}^{N}
       \frac{|u_j(v)|^2}{(x-\lambda_j)^2}<0.
\]
At \(x=\lambda_j\), the polynomial on the left equals
\(-|u_j(v)|^2p'(\lambda_j)^2<0\), so it has no real zero.
For a real symmetric matrix, \(e_v\) is cyclic if and only if
the spectrum is simple and \(e_v\) has a nonzero projection onto
each eigenspace. Hence coprimality implies cyclicity. Conversely,
if \(e_v\) is cyclic, then \(u_j(v)\neq0\) and \(p'(\lambda_j)\neq0\)
for every \(j\). The identity for \(q(\lambda_j)\) therefore gives
\(\gcd(p,q)=1\).
\end{proof}

\begin{lemma}
\label{lem:resultant-walk-identity}
Let \(M\) be a real symmetric matrix of order \(N\), and let
\(v\in\{1,\ldots,N\}\). Set
\[
  p(x)=\phi(M;x),\qquad q(x)=\phi(M_{\hat v};x).
\]
Then
\begin{equation}\label{eq:principal-minor-krylov}
  \det W(M,e_v)^2
  =(-1)^{\binom N2}\Res(p,q).
\end{equation}
\end{lemma}

\begin{proof}
We first assume that \(M\) has distinct eigenvalues
\(\lambda_1,\ldots,\lambda_N\), with an orthonormal eigenbasis
\(u_1,\ldots,u_N\). Set \(c_j=u_j^Te_v\).
Spectral decomposition gives
\[
  W(M,e_v)
  =[u_1,\ldots,u_N]\diag(c_1,\ldots,c_N)
    (\lambda_j^{k-1})_{j,k=1}^{N}.
\]
The last factor is the Vandermonde matrix associated with the eigenvalues.
The spectral expansion \eqref{eq:principal-resolvent} gives
\(q(\lambda_j)=c_j^2p'(\lambda_j)\) as in Lemma~\ref{lem:wronskian-resultant-cyclic}.
The eigenvector matrix is orthogonal, so the above
factorization gives
\[
  \det W(M,e_v)^2
  =\left(\prod_{j=1}^{N}c_j^2\right)
   \prod_{1\leq i<j\leq N}(\lambda_j-\lambda_i)^2.
\]
Since \(p\) is monic, the product formula for the resultant gives
\[
  \Res(p,q)
  =\prod_{j=1}^{N}q(\lambda_j)
  =\left(\prod_{j=1}^{N}c_j^2\right)
   \prod_{j=1}^{N}p'(\lambda_j).
\]
Using
\[
  \prod_{j=1}^{N}p'(\lambda_j)
  =(-1)^{\binom N2}
   \prod_{1\leq i<j\leq N}(\lambda_j-\lambda_i)^2,
\]
we obtain \eqref{eq:principal-minor-krylov}.

Finally, matrices with simple spectrum are dense in the space of
real symmetric matrices, and both sides of
\eqref{eq:principal-minor-krylov} are polynomials in the matrix
entries. Continuity therefore proves the identity for every real
symmetric \(M\).
\end{proof}

For a polynomial \(h(t)=\sum_k c_kt^k\) and a square matrix \(S\),
we use \(h(S)=\sum_k c_kS^k\), with \(S^0=I\).
The rooted-product criterion follows from these identities and the
following factorization.

\begin{theorem}[Wang, Shen, and Mao
  {\cite[Theorem~3]{WangShenMao2026}}]
\label{thm:wang-shen-mao-rooted}
Let \(M\) be \(A\) or \(Q\). Let \(G\) have
order \(n\), and let \(H^{(v)}\) be a rooted graph of order \(m\). Set
\[
  B_M(\lambda)=M(H)+\lambda E_{vv}.
\]
Then
\begin{equation}\label{eq:wang-shen-mao-rooted}
\begin{split}
  \det W_M(G\circ H^{(v)})
  ={}&\pm
  \Res\bigl(
    \phi(M(H);x),
    \phi(M(H)_{\hat v};x)
  \bigr)^{n(n-1)/2}\\
  &{}\cdot\det h_{M,H,v}(M(G))
  \cdot\bigl(\det W_M(G)\bigr)^m,
\end{split}
\end{equation}
where
\[
  h_{M,H,v}(\lambda)
  =\det W(B_M(\lambda))
  =\det\bigl[
    \boldsymbol{1},B_M(\lambda)\boldsymbol{1},\ldots,
    B_M(\lambda)^{m-1}\boldsymbol{1}
  \bigr].
\]
\end{theorem}

The proof of the cited factorization uses the Kronecker decomposition
and the rank-one perturbations \(M(H)+\lambda E_{vv}\).
Equation~\eqref{eq:Aalpha-rooted} gives this decomposition for \(A_{\alpha}\),
up to a simultaneous permutation of rows and columns, so the argument
also yields \eqref{eq:wang-shen-mao-rooted} for \(M=A_{\alpha}\).
Define the resultant associated with the root by
\[
  R_{\alpha}(H,v)
  =
  \Res\bigl(
    \phi(A_{\alpha}(H);x),
    \phi\bigl(\bigl(A_{\alpha}(H)\bigr)_{\hat v};x\bigr)
  \bigr).
\]
The corresponding walk-determinant polynomial is
\[
  h_{\alpha,H,v}(t)
  =
  \det W\bigl(A_{\alpha}(H)+tE_{vv}\bigr).
\]

The three conditions in the following corollary are those of
Shan and Liu~\cite[Theorem~6.1]{ShanLiu2025}, written using the
walk-determinant polynomial.

\begin{corollary}
\label{cor:rooted-general}
Let \(G\) have order \(n\geq2\), and let \(H^{(v)}\) have order
\(m\geq2\).
Then \(G\circ H^{(v)}\) is \(A_{\alpha}\)-controllable if and
only if the following three conditions hold:
\begin{enumerate}
  \item \(G\) is \(A_{\alpha}\)-controllable;
  \item \(v\) is an \(A_{\alpha}\)-Wronskian vertex of \(H\);
  \item \(h_{\alpha,H,v}(\lambda)\neq0\) for every
        \(\lambda\in\Spec(A_{\alpha}(G))\).
\end{enumerate}
If the second condition holds, then \(h_{\alpha,H,v}\not\equiv0\).
Thus the third condition excludes only finitely many real numbers.
\end{corollary}

\begin{proof}
Apply the \(A_{\alpha}\) form of \eqref{eq:wang-shen-mao-rooted},
as justified immediately after Theorem~\ref{thm:wang-shen-mao-rooted}.
Since \(n\geq2\), the exponent \(n(n-1)/2\) is positive, so the
resultant factor is nonzero if and only if the resultant itself is nonzero.
Lemma~\ref{lem:wronskian-resultant-cyclic} identifies the nonzero
resultant condition with the root condition.
Let \(\lambda_1,\ldots,\lambda_n\) be the eigenvalues of
\(A_{\alpha}(G)\), counted with multiplicity. Then
\[
  \det h_{\alpha,H,v}(A_{\alpha}(G))
  =
  \prod_{j=1}^{n}h_{\alpha,H,v}(\lambda_j).
\]
Substituting this product into the determinant factorization
shows that the determinant is nonzero exactly when the three conditions hold.

It remains to show that \(h_{\alpha,H,v}\not\equiv0\) when
the second condition holds. Set \(M=A_{\alpha}(H)\),
\(p(x)=\phi(M;x)\), and \(q(x)=\phi(M_{\hat v};x)\).
Then \(\gcd(p,q)=1\). Define
\[
  B(t)=M+tE_{vv},\qquad
  \rho(x)=\boldsymbol1^T\operatorname{adj}(xI-M)e_v.
\]
Here \(\operatorname{adj}(S)\) is the adjugate matrix, the transpose
of the cofactor matrix, satisfying \(S\operatorname{adj}(S)=\det(S)I\).
The polynomial \(\rho\) is monic of degree \(m-1\), since
the leading coefficient of \(\operatorname{adj}(xI-M)\) is \(I\)
and \(\boldsymbol1^Te_v=1\).
Moreover,
\[
  \phi(B(t);x)=p(x)-tq(x).
\]
Since \(\gcd(p-tq,q)=1\), Lemma~\ref{lem:wronskian-resultant-cyclic}
shows that \(e_v\) is cyclic for \(B(t)\).
Thus \(B(t)\) has simple spectrum, and every eigenvector
has a nonzero \(v\)-coordinate.

Suppose that \(h_{\alpha,H,v}\equiv0\).
For each \(t>0\), noncyclicity of \(\boldsymbol1\) for \(B(t)\)
gives an eigenvalue \(\mu_t\) and an eigenvector \(\boldsymbol y_t\)
orthogonal to \(\boldsymbol1\), with \(y_t(v)\neq0\).
The eigenvector equation gives
\[
  (\mu_t I-M)\boldsymbol y_t=t\,y_t(v)e_v.
\]
Multiplying by \(\operatorname{adj}(\mu_t I-M)\) and then by
\(\boldsymbol1^T\) yields
\[
  0=p(\mu_t)\boldsymbol1^T\boldsymbol y_t
   =t\,y_t(v)\rho(\mu_t).
\]
Hence \(\rho(\mu_t)=0\).
If \(t\neq s\) and \(\mu_t=\mu_s=\mu\), then
\(p(\mu)=tq(\mu)=sq(\mu)\), so \(p(\mu)=q(\mu)=0\),
contrary to \(\gcd(p,q)=1\).
Thus \(\rho\) would have infinitely many distinct zeros,
contradicting its degree. Therefore \(h_{\alpha,H,v}\not\equiv0\).
\end{proof}

\subsection{Even-order paths rooted at a central vertex}

For either central vertex of \(P_{2r}\), the two sides have orders
\(r-1\) and \(r\). Theorem~\ref{thm:main} therefore shows that
both vertices are \(A_{\alpha}\)-Wronskian for every
\(0\leq\alpha<1\). Roots at end vertices also have this property and
were used in~\cite[Theorem~4.2]{TianWuCui2024}.
Here the consecutive block orders allow us to evaluate both
path-dependent factors explicitly.

Write \(P_{2r}^{(j)}\) for \(P_{2r}\) rooted at vertex \(j\), and set
\[
  a_r=\left\lfloor\frac r2\right\rfloor,
  \qquad
  b_r=\left\lfloor\frac{r-1}{2}\right\rfloor.
\]

\begin{theorem}\label{thm:even-center}
Let \(r\geq1\) and \(0\leq\alpha<1\), and put \(\beta=1-\alpha\).
For either central vertex \(v\in\{r,r+1\}\),
\begin{equation}\label{eq:even-center-resultant}
  R_{\alpha}(P_{2r},v)
  =
  (-1)^r\beta^{\,2r(2r-1)}.
\end{equation}
At root \(r\),
\begin{equation}\label{eq:even-center-h}
\begin{split}
  h_{\alpha,P_{2r},r}(t)
  ={}&
  \pm\beta^{\,2(r-1)^2}t^r\\
  &{}\cdot
  (\alpha t+2\alpha-1)^{a_r}
  (\alpha t+1-2\alpha)^{b_r}.
\end{split}
\end{equation}
The sign in \eqref{eq:even-center-h} depends only on \(r\).
In the natural vertex order,
\(h_{\alpha,P_{2r},r+1}(t)=(-1)^r h_{\alpha,P_{2r},r}(t)\).
\end{theorem}

At \(\alpha=0\), substituting these formulas into
\eqref{eq:wang-shen-mao-rooted} recovers the central-root case of
Yan and Wang~\cite[Theorem~3]{YanWang2026}.

Substituting Theorem~\ref{thm:even-center} into the rooted-product
factorization gives, for a graph \(G\) of order \(n\) and either
\(v\in\{r,r+1\}\),
\begin{equation}\label{eq:even-center-product}
\begin{split}
 &\det\Walpha(G\circ P_{2r}^{(v)})\\
 &\quad={}
 \pm\beta^{E_{n,r}}
 \det(A_{\alpha}(G))^r\\
 &\qquad{}\cdot
 \det\bigl(\alpha A_{\alpha}(G)+(2\alpha-1)I_n\bigr)^{a_r}\\
 &\qquad{}\cdot
 \det\bigl(\alpha A_{\alpha}(G)+(1-2\alpha)I_n\bigr)^{b_r}
 \bigl(\det\Walpha(G)\bigr)^{2r},
\end{split}
\end{equation}
where
\[
  E_{n,r}
  =r(2r-1)n(n-1)+2(r-1)^2n.
\]

\begin{corollary}[Central-root criterion]
\label{cor:even-center-controllable}
Let \(G\) be a graph, \(r\geq1\), \(0\leq\alpha<1\), and
\(v\in\{r,r+1\}\). Then \(G\circ P_{2r}^{(v)}\) is
\(A_{\alpha}\)-controllable if and only if all the following conditions hold:
\begin{enumerate}[label=\textup{(\roman*)}]
  \item \(G\) is \(A_{\alpha}\)-controllable and \(A_{\alpha}(G)\)
  is nonsingular;
  \item if \(\alpha>0\) and \(r\geq2\), then
  \((1-2\alpha)/\alpha\) is not an eigenvalue of \(A_{\alpha}(G)\);
  \item if \(\alpha>0\) and \(r\geq3\), then
  \((2\alpha-1)/\alpha\) is not an eigenvalue of \(A_{\alpha}(G)\).
\end{enumerate}
For \(r=1\) or \(\alpha\in\{0,1/2\}\), condition \textup{(i)} alone suffices.
\end{corollary}

\begin{proof}
If \(G\) has one vertex, then \(A_{\alpha}(G)=[0]\), so condition
\textup{(i)} fails. The product is \(P_{2r}\). Path reversal fixes
\(\boldsymbol1\) and commutes with \(A_{\alpha}(P_{2r})\), so the walk
matrix has identical rows at reflected vertices and is singular.

Now suppose that \(G\) has at least two vertices.
Since \(\beta>0\), formula \eqref{eq:even-center-product}
is nonzero exactly when \(G\) is \(A_{\alpha}\)-controllable,
\(A_{\alpha}(G)\) is nonsingular, and each factor with a positive
exponent is nonzero. We have \(a_r>0\) if and only if \(r\geq2\),
and \(b_r>0\) if and only if \(r\geq3\).
For \(\alpha>0\), the two factors exclude
\((1-2\alpha)/\alpha\) and \((2\alpha-1)/\alpha\), respectively.
At \(\alpha=0\), they are nonzero constants; at \(\alpha=1/2\),
both excluded values are zero. Thus these conditions are either
absent or already included in condition \textup{(i)}.
\end{proof}

Suppose \(1/2<\alpha<1\) and \(G\) has minimum degree at least \(2\).
Since \(Q(G)\) is positive semidefinite, the identity
\(A_{\alpha}(G)=(1-\alpha)Q(G)+(2\alpha-1)D(G)\) yields
\(\lambda_{\min}(A_{\alpha}(G))\geq2(2\alpha-1)>(2\alpha-1)/\alpha\).
The spectral exclusions are then automatic: for every \(r\geq1\) and
either central root, the product is \(A_{\alpha}\)-controllable if and
only if the base graph is.

\begin{remark}
At \(\alpha=1/2\), both additional excluded values are zero.
Since \(A_{1/2}(G)=Q(G)/2\),
Corollary~\ref{cor:even-center-controllable} becomes:
\(G\circ P_{2r}^{(v)}\) is \(Q\)-controllable if and only if
\(G\) is \(Q\)-controllable and \(Q(G)\) is nonsingular.
For connected \(G\), nonsingularity of \(Q(G)\) is equivalent to
\(G\) being non-bipartite.

Scaling the walk columns and applying \eqref{eq:even-center-h} gives,
for either \(v\in\{r,r+1\}\),
\[
  \det W(Q(P_{2r})+tE_{vv})
  =2^{\binom{2r}{2}}h_{1/2,P_{2r},v}(t/2)
  =\pm t^{2r-1}.
\]
\end{remark}

We now prove Theorem~\ref{thm:even-center} using two resultant calculations.
For the walk-determinant formula, we first work with \(0<\alpha<1\)
and then recover \(\alpha=0\) by a polynomial identity.
We suppress the variable \(x\).
Set
\[
  D_r=f_{r-1}f_r,
  \qquad
  \Phi_{2r}=f_r^2-\beta^2f_{r-1}^2.
\]
Thus \(D_r\) is the characteristic polynomial of the root-deleted
principal submatrix and
\(\Phi_{2r}=\phi(A_{\alpha}(P_{2r});x)\).
For the latter identity, split the matrix at the central edge,
keep its diagonal entries unchanged, and reverse one half.
The two \(B_r\) blocks are coupled by \(\beta E_{rr}\).
The subspaces of vectors \((u,u)\) and \((u,-u)\) give the
characteristic factors \(f_r-\beta f_{r-1}\) and
\(f_r+\beta f_{r-1}\), respectively, proving the formula for \(\Phi_{2r}\).
Reducing \(\Phi_{2r}\) modulo the two factors of \(D_r\) gives
\[
  \Phi_{2r}\equiv f_r^2\pmod{f_{r-1}},
  \qquad
  \Phi_{2r}\equiv-\beta^2f_{r-1}^2\pmod{f_r}.
\]
To express eigenvector coordinate sums on the two sides of the root,
for \(j\geq0\) define \(\eta_j=f_j+\sigma_j\), where
\[
  \sigma_0=0,
  \qquad
  \sigma_j=\sum_{k=0}^{j-1}\beta^{j-k}f_k\quad(j\geq1).
\]
Then \(\sigma_j=\beta\eta_{j-1}\) for \(j\geq1\). We use the
factorization
\begin{equation}\label{eq:even-center-Psi}
  \Psi_r
  :=f_{r-1}f_r+\sigma_{r-1}f_r+\sigma_rf_{r-1}
  =\eta_{r-1}(f_r+\beta f_{r-1}).
\end{equation}

The first lemma expresses the walk-determinant polynomial as a resultant
involving \(\Psi_r\). By \eqref{eq:even-center-Psi}, it then
suffices to evaluate the contributions of \(\eta_{r-1}\) and
\(f_r+\beta f_{r-1}\). The latter gives \(t^r\); the second lemma
gives the two remaining factors in \eqref{eq:even-center-h}.

\begin{lemma}\label{lem:center-coordinate-sum}
Assume \(0<\alpha<1\), let \(t\in\R\), and put
\(C_r(t)=A_{\alpha}(P_{2r})+tE_{rr}\). If \(\mu\) is an eigenvalue of
\(C_r(t)\), then every corresponding eigenvector \(\boldsymbol{y}\)
satisfies \(y_r\neq0\). After normalizing \(\boldsymbol{y}\) so that
\(y_r=1\), we have
\[
  \boldsymbol{1}^T\boldsymbol{y}=\frac{\Psi_r(\mu)}{D_r(\mu)}.
\]
Moreover,
\begin{equation}\label{eq:even-center-h-resultant}
  h_{\alpha,P_{2r},r}(t)
  =\pm |R_{\alpha}(P_{2r},r)|^{-1/2}
  \Res\bigl(\Psi_r,\Phi_{2r}-tD_r\bigr).
\end{equation}
\end{lemma}

\begin{proof}
Set \(C=C_r(t)\) and recall \(D_r=f_{r-1}f_r\).
The characteristic polynomial of \(C\) is
\[
  p_t=\Phi_{2r}-tD_r.
\]
Since \(\beta>0\), the recurrence \eqref{eq:frecurrence} gives
\(\gcd(f_{r-1},f_r)=1\). Since \(p_t-\Phi_{2r}=-tD_r\), the
congruences for \(\Phi_{2r}\) also hold for \(p_t\), giving
\(\gcd(p_t,D_r)=1\). By Lemma~\ref{lem:wronskian-resultant-cyclic},
\(e_r\) is cyclic for \(C\). In particular, the spectrum is simple,
every eigenvector is nonzero at the root, and \(D_r(\mu)\neq0\)
at every eigenvalue \(\mu\).

Normalize a \(\mu\)-eigenvector \(\boldsymbol y\) by \(y_r=1\).
Induction along the two sides of the path, using the eigenvector
equations and \eqref{eq:frecurrence}, gives
\[
  y_j=\frac{\beta^{r-j}f_{j-1}(\mu)}{f_{r-1}(\mu)}
  \quad(1\leq j\leq r-1),\qquad
  y_{r+q}=\frac{\beta^qf_{r-q}(\mu)}{f_r(\mu)}
  \quad(1\leq q\leq r).
\]
Summing and using the definitions of \(\sigma_j\) and \(\Psi_r\)
gives
\[
  \boldsymbol1^T\boldsymbol y
  =1+\frac{\sigma_{r-1}(\mu)}{f_{r-1}(\mu)}
     +\frac{\sigma_r(\mu)}{f_r(\mu)}
  =\frac{\Psi_r(\mu)}{D_r(\mu)}.
\]
This proves the first assertion.

To derive the determinant identity, compare \(W(C,\boldsymbol1)\)
and \(W(C,e_r)\). Let \(\mu_1,\ldots,\mu_{2r}\) be the eigenvalues
of \(C\), with an orthonormal eigenbasis \(u_1,\ldots,u_{2r}\).
Applying the coordinate-sum formula to \(u_j/u_j(r)\) gives
\[
  \frac{\boldsymbol1^Tu_j}{u_j(r)}
  =\frac{\Psi_r(\mu_j)}{D_r(\mu_j)}.
\]
The two spectral factorizations have the same eigenvector and
Vandermonde matrices, so
\[
  \frac{h_{\alpha,P_{2r},r}(t)}{\det W(C,e_r)}
  =\prod_{j=1}^{2r}\frac{\boldsymbol1^Tu_j}{u_j(r)}
  =\frac{\Res(p_t,\Psi_r)}{\Res(p_t,D_r)},
\]
where \(p_t\) is monic and the denominator is nonzero because
\(e_r\) is cyclic. Equation~\eqref{eq:principal-minor-krylov} gives
\[
  |\det W(C,e_r)|^2=|\Res(p_t,D_r)|,
  \qquad
  |h_{\alpha,P_{2r},r}(t)|
  =\frac{|\Res(p_t,\Psi_r)|}{|\Res(p_t,D_r)|^{1/2}}.
\]
Finally, \(p_t=\Phi_{2r}\) at every zero of \(D_r\), so
\[
  \Res(p_t,D_r)=\Res(\Phi_{2r},D_r)=R_{\alpha}(P_{2r},r).
\]
Since \(\deg p_t=2r\) is even, interchanging the two arguments
of \(\Res(p_t,\Psi_r)\) leaves it unchanged.
This proves \eqref{eq:even-center-h-resultant}.
\end{proof}

\begin{lemma}\label{lem:eta-resultant}
Assume \(0<\alpha<1\) and let \(t\in\R\). Then
\[
\begin{split}
 &\left|\Res
  \bigl(\eta_{r-1},\Phi_{2r}-tD_r\bigr)\right|\\
 &\qquad={}
 \beta^{2(r-1)^2}
 |\alpha t+2\alpha-1|^{a_r}
 |\alpha t+1-2\alpha|^{b_r}.
\end{split}
\]
\end{lemma}

\begin{proof}
For \(r=1\), the formula follows from \(\eta_0=1\).
Assume \(r\geq2\) and put \(m=r-1\). We use the product formula
for the resultant: first compute the product of \(D_r\) over the zeros
of \(\eta_m\), then split \(\eta_m\) into factors on whose zeros
\(\Phi_{2r}/D_r\) is constant. These constants will determine the
two linear factors in \(t\).

First, we compute the required product. More generally, for every \(j\geq0\),
\begin{equation}\label{eq:eta-resultant}
  \left|\Res(\eta_j,f_jf_{j+1})\right|
  =\alpha^j\beta^{2j^2}.
\end{equation}
The case \(j=0\) follows from \(\eta_0=1\). Set \(x=\beta z+2\alpha\).
In calculations in \(z\), all polynomials in \(x\) are evaluated at
this expression. Define
\[
  S_j(z)=\sum_{k=0}^{j}Y_k(z).
\]
Then \(f_j(x)=\beta^jY_j(z)\) and
\(\eta_j(x)=\beta^jS_j(z)\).
The recurrence for \(Y_j\) gives
\[
  S_0=1,\qquad S_1=z+1+\tau,\qquad
  S_j=zS_{j-1}-S_{j-2}+1+\tau\quad(j\geq2).
\]
Consequently, for every \(j\geq1\),
\[
  S_j^2+S_{j-1}^2-zS_jS_{j-1}
  -(1+\tau)(S_j+S_{j-1})=-\tau.
\]
Indeed, substitution of the recurrence shows that the left-hand
side is unchanged when \(j\) is replaced by \(j-1\), for \(j\geq2\),
and its value at \(j=1\) is \(-\tau\).
Since
\[
  Y_j=S_j-S_{j-1},\qquad
  Y_{j+1}=(z-1)S_j-S_{j-1}+1+\tau,
\]
reduction modulo \(S_j\) gives
\[
  Y_jY_{j+1}
  \equiv S_{j-1}(S_{j-1}-1-\tau)
  \equiv-\tau\pmod{S_j}.
\]
Returning to the variable \(x\) and using \(\beta\tau=\alpha\),
we obtain
\[
  f_jf_{j+1}\equiv-\alpha\beta^{2j}\pmod{\eta_j}.
\]
The polynomial \(\eta_j\) is monic of degree \(j\). The product
formula for the resultant, with zeros counted with multiplicity,
therefore gives
\[
  \Res(\eta_j,f_jf_{j+1})=(-\alpha\beta^{2j})^j.
\]
Taking absolute values proves \eqref{eq:eta-resultant}.
For \(j=m\), this is the required product because \(D_r=f_mf_{m+1}\).
No assumption that the zeros of \(\eta_j\) are simple is needed.

Next, we seek a factorization of \(\eta_m\) for which
\(\Phi_{2r}/D_r\) is constant on the zeros of each factor.
To factor \(S_m(z)\), we use \(z=4c^2-2\).
For \(z=2\cos\theta\), this is the half-angle substitution
\(c=\cos(\theta/2)\). We regard this as a polynomial substitution
with \(c\) a formal variable, so it imposes no restriction on \(z\).
For \(k\geq0\), the double-angle identity
\[
  2c\,U_k(2c^2-1)
  =U_k(c)\bigl(U_{k+1}(c)-U_{k-1}(c)\bigr)
\]
follows from the sine representation of \(U_k\) with
\(c=\cos\varphi\), \(0<\varphi<\pi/2\), and hence holds as a
polynomial identity.
The finite-sum identity
\[
 \sum_{k=0}^{m}
 \left[U_k(2c^2-1)+\tau U_{k-1}(2c^2-1)\right]
 =\frac{U_m(c)\bigl(U_{m+1}(c)+\tau U_{m-1}(c)\bigr)}{2c}
\]
can be verified by telescoping. Set
\(N_k=U_k(c)\bigl(U_{k+1}(c)+\tau U_{k-1}(c)\bigr)\).
Then \(N_0=2c\), and the double-angle identity gives, for \(k\geq1\),
\[
\begin{aligned}
 N_k-N_{k-1}
 &=U_k(c)\bigl(U_{k+1}(c)-U_{k-1}(c)\bigr)\\
 &\quad+\tau U_{k-1}(c)\bigl(U_k(c)-U_{k-2}(c)\bigr)\\
 &=2c\left[U_k(2c^2-1)+\tau U_{k-1}(2c^2-1)\right].
\end{aligned}
\]
Summing over \(1\leq k\leq m\) proves the identity.
We assign the division by \(2c\) to the odd numerator factor and define
\[
 (\mathcal P_m(z),\mathcal Q_m(z))=
 \begin{cases}
 \left(
   U_m(c),
   \dfrac{U_{m+1}(c)+\tau U_{m-1}(c)}{2c}
 \right),&m\text{ even},\\[4mm]
 \left(
   \dfrac{U_m(c)}{2c},
   U_{m+1}(c)+\tau U_{m-1}(c)
 \right),&m\text{ odd}.
 \end{cases}
\]
Since \(U_j(-c)=(-1)^jU_j(c)\), every numerator divided by \(2c\)
is odd. Each division is therefore exact:
after cancellation they are even polynomials in \(c\), defined
also at \(c=0\). Thus both factors are polynomials in
\(z=4c^2-2\). Since \(U_j(c)\) has leading term \((2c)^j\),
they are monic in \(z\), of degrees \(b_r\) and \(a_r\), respectively.
Substituting them into the finite-sum identity gives
\begin{equation}\label{eq:eta-polynomial-split}
  \beta^{-m}\eta_m(\beta z+2\alpha)=\mathcal P_m(z)\mathcal Q_m(z).
\end{equation}
To find the two constant values, write \(D_r\) and \(\Phi_{2r}\)
in terms of \(Y_m\) and \(Y_{m+1}\):
\[
  D_r=\beta^{2m+1}Y_mY_{m+1},
  \qquad
  \Phi_{2r}
  =\beta^{2m+2}\bigl(Y_{m+1}^2-Y_m^2\bigr).
\]
Using \(\alpha=\beta\tau\) and \(\beta(1+\tau)=1\), we obtain
\[
\begin{split}
  \alpha\Phi_{2r}-(2\alpha-1)D_r
  &=
  \beta^{2m+3}
  (Y_m+\tau Y_{m+1})
  (Y_{m+1}-\tau Y_m),\\
  \alpha\Phi_{2r}+(2\alpha-1)D_r
  &=
  \beta^{2m+3}
  (Y_{m+1}+\tau Y_m)
  (\tau Y_{m+1}-Y_m).
\end{split}
\]

The substitution \(z/2=2c^2-1\) and the recurrence
\(U_{j+1}(c)=2cU_j(c)-U_{j-1}(c)\) give
\[
  Y_m+\tau Y_{m+1}
  =
  \bigl(1+\tau^2+2\tau(2c^2-1)\bigr)
  U_m(2c^2-1).
\]
The double-angle identity, the definition of \(\mathcal P_m\),
and the parity of \(U_m\) now give
\[
  \mathcal P_m(z)\mid Y_m(z)+\tau Y_{m+1}(z).
\]

For the other factor, the same recurrence gives
\[
\begin{split}
  Y_{m+1}+\tau Y_m
  ={}&
  \frac{U_{m+1}(c)+\tau U_{m-1}(c)}{2c}
  \Bigl(
    U_{m+2}(c)-U_m(c)
    +\tau\bigl(U_m(c)-U_{m-2}(c)\bigr)
  \Bigr),
\end{split}
\]
so the definition of \(\mathcal Q_m\) gives
\[
  \mathcal Q_m(z)\mid Y_{m+1}(z)+\tau Y_m(z).
\]
Consequently,
\begin{equation}\label{eq:eta-split-congruences}
\begin{split}
  \alpha\Phi_{2r}-(2\alpha-1)D_r
    &\equiv0\pmod{\mathcal P_m(z)},\\
  \alpha\Phi_{2r}+(2\alpha-1)D_r
    &\equiv0\pmod{\mathcal Q_m(z)}.
\end{split}
\end{equation}
By \eqref{eq:eta-resultant} with \(j=m\),
\(\Res(\eta_m,D_r)\neq0\). Thus \(D_r\) is nonzero at every
zero of either factor, and the congruences give
\[
  \frac{\Phi_{2r}}{D_r}=
  \begin{cases}
    (2\alpha-1)/\alpha, & \text{at the zeros of }\mathcal P_m,\\[2pt]
    (1-2\alpha)/\alpha, & \text{at the zeros of }\mathcal Q_m.
  \end{cases}
\]
Each zero is counted with its multiplicity in each factor,
including zeros common to both factors. All preceding calculations
are polynomial identities, so repeated zeros and \(c=0\) cause no exception.

Finally, evaluate the product for \(\Phi_{2r}-tD_r\).
At each zero, its value is \(D_r\) times the following multiplier:
\[
  \frac{\Phi_{2r}-tD_r}{D_r}
  =-\frac1\alpha
  \begin{cases}
    \alpha t+1-2\alpha, & \text{at the zeros of }\mathcal P_m,\\[2pt]
    \alpha t+2\alpha-1, & \text{at the zeros of }\mathcal Q_m.
  \end{cases}
\]
Thus \(\mathcal P_m\) and \(\mathcal Q_m\) contribute
\(|\alpha t+1-2\alpha|^{b_r}\) and
\(|\alpha t+2\alpha-1|^{a_r}\), respectively, together with \(\alpha^{-m}\).
By \eqref{eq:eta-resultant}, the product of the \(D_r\)-values
has absolute value \(\alpha^m\beta^{2m^2}\).
The powers of \(\alpha\) cancel, proving the formula since
\(a_r+b_r=m=r-1\).
\end{proof}

\begin{proof}[Proof of Theorem~\ref{thm:even-center}]
We first compute \(R_{\alpha}(P_{2r},r)\).
Recall \(D_r=f_{r-1}f_r\). Put
\(T_j=|\Res(f_j,f_{j-1})|\).
Since \(f_{j-1}\) is monic of degree \(j-1\), reduction modulo
\(f_{j-1}\) in \eqref{eq:frecurrence} gives
\[
  T_1=1,\qquad T_j=\beta^{2(j-1)}T_{j-1}\quad(j\geq2).
\]
Hence
\[
  |\Res(f_j,f_{j-1})|=\beta^{j(j-1)}.
\]
Since \(\deg\Phi_{2r}=2r\) is even, we may interchange the arguments
of \(R_{\alpha}(P_{2r},r)=\Res(\Phi_{2r},D_r)\).
By multiplicativity and the congruences for \(\Phi_{2r}\), the factors
of \(\Res(D_r,\Phi_{2r})\) are \(T_r^2\) and \((-\beta^2)^rT_r^2\).
It follows that
\[
  R_{\alpha}(P_{2r},r)
  =(-\beta^2)^rT_r^4
  =(-1)^r\beta^{2r(2r-1)}.
\]

To compute \(h_{\alpha,P_{2r},r}\), first assume \(0<\alpha<1\).
Set \(\Lambda_r=f_r+\beta f_{r-1}\), the second factor in
\eqref{eq:even-center-Psi}. Modulo \(\Lambda_r\),
\(\Phi_{2r}=0\) and \(-tD_r=\beta t f_{r-1}^2\), so
\[
\begin{split}
 \left|\Res(\Lambda_r,\Phi_{2r}-tD_r)\right|
  &=\beta^r|t|^r
   \left|\Res(f_r,f_{r-1})\right|^2\\
 &=\beta^{r(2r-1)}|t|^r.
\end{split}
\]
By multiplicativity, this identity and
Lemma~\ref{lem:eta-resultant} give the resultant involving
\(\Psi_r\). Substituting the result into
\eqref{eq:even-center-h-resultant} gives
\[
\begin{split}
 |h_{\alpha,P_{2r},r}(t)|
 ={}&\beta^{2(r-1)^2}|t|^r
 |\alpha t+2\alpha-1|^{a_r}\\
 &{}\cdot|\alpha t+1-2\alpha|^{b_r}.
\end{split}
\]
To remove the absolute values, define the polynomial
\[
  H_r(\alpha,t)=
  \beta^{2(r-1)^2}t^r
  (\alpha t+2\alpha-1)^{a_r}
  (\alpha t+1-2\alpha)^{b_r}.
\]
The preceding identity gives
\(h_{\alpha,P_{2r},r}(t)^2=H_r(\alpha,t)^2\) for
\(0<\alpha<1\) and every real \(t\).
Both sides are polynomials in \(\alpha\) and \(t\), so the
equality is a polynomial identity and also holds at \(\alpha=0\).
Since \(\R[\alpha,t]\) is an integral domain,
\(h_{\alpha,P_{2r},r}=\pm H_r\), with a sign independent of
\(\alpha\) and \(t\). This proves \eqref{eq:even-center-h}.
Path reversal preserves the two characteristic polynomials in the
resultant. It permutes the rows of the walk matrix by \(r\) transpositions,
giving the stated sign relation between the two roots.
\end{proof}

\subsection{Uniform base graphs for even-order paths}

By Corollary~\ref{cor:even-center-controllable}, the test at \(r=3\)
includes every spectral exclusion that can occur, since \(a_3=b_3=1\).
For a fixed \(\alpha\), it thus suffices to test \(G\circ P_6^{(3)}\)
to cover all \(r\geq1\) and both central roots.

Writing \(\alpha=1/x\), the two additional excluded values become
\(x-2\) and \(2-x\). A row-sum bound excludes both for large \(x\),
giving a sufficient condition valid for every \(r\geq1\).

\noindent\begin{minipage}{\linewidth}
\begin{theorem}
\label{thm:eventual-reciprocal}
Let \(G\) be a graph of order \(n\), and set
\[
  M(x)=D(G)+(x-1)A(G),\qquad
  F_G(x)=\det W(M(x)),\qquad g_G(x)=\det M(x).
\]
Let \(\Delta(G)\) denote the maximum degree of \(G\). Suppose that
\(F_Gg_G\not\equiv0\). Then, for every real number
\(x>\Delta(G)+2\) satisfying \(F_G(x)g_G(x)\neq0\), both
\(G\circ P_{2r}^{(r)}\) and \(G\circ P_{2r}^{(r+1)}\)
are \(A_{1/x}\)-controllable for all \(r\geq1\).

Consequently, there exists \(\varepsilon>0\) such that both products
are \(A_{\alpha}\)-controllable for every
\(0<\alpha<\varepsilon\) and every \(r\geq1\).
\end{theorem}
\end{minipage}

\begin{proof}
Set \(N=\binom{n}{2}\). For \(x\) satisfying the hypotheses,
\[
  A_{1/x}(G)=x^{-1}M(x),
  \qquad
  \det W(A_{1/x}(G))=x^{-N}\det W(M(x))
\]
and \(F_G(x)g_G(x)\neq0\) show that \(G\) is
\(A_{1/x}\)-controllable and \(A_{1/x}(G)\) is nonsingular.

Let \(\Delta=\Delta(G)\). The row sum of the nonnegative matrix
\(A_{1/x}(G)\) at \(v\) is \(d(v)\), so every eigenvalue satisfies
\(|\lambda|\leq\Delta\). Since \(x>\Delta+2\), neither \(x-2\)
nor \(2-x\) is an eigenvalue.
These are the additional values in
Corollary~\ref{cor:even-center-controllable} at \(\alpha=1/x\).
Together with controllability and nonsingularity, their exclusion
proves the assertion for every \(r\geq1\).

Finally, \(F_Gg_G\) is a nonzero polynomial, so it has finitely
many real zeros. Choose \(X>\Delta+2\) larger than all its positive
zeros, and set \(\varepsilon=1/X\).
For \(0<\alpha<\varepsilon\), we have \(x=1/\alpha>X\), and the
first assertion applies.
\end{proof}

We next cover the full parameter range using a graph and its complement.
Let \(G_*\) be the graph on
\(\{1,2,\ldots,6\}\) with
\begin{equation}\label{eq:seed-edges}
  E(G_*)=
  \bigl\{
    \{1,3\},\{1,4\},\{1,6\},
    \{2,3\},\{2,5\},\{3,4\}
  \bigr\}.
\end{equation}

\begin{theorem}[A six-vertex graph and its complement]
\label{thm:two-seed-cover}
Let \(G_*\) be the graph in \eqref{eq:seed-edges}.
There is a finite set \(\mathcal E\subset[0,1)\), with
\(\mathcal E\cap[0,1/5)=\varnothing\), such that:
\begin{enumerate}[label=\textup{(\roman*)}]
  \item for each \(\alpha\in[0,1)\setminus\mathcal E\),
  \(G_*\circ P_{2r}^{(v)}\) is \(A_{\alpha}\)-controllable for
  every \(r\geq1\) and both \(v\in\{r,r+1\}\);
  \item for each \(\alpha\in\mathcal E\),
  \(\overline{G_*}\circ P_{2r}^{(v)}\) is
  \(A_{\alpha}\)-controllable for every \(r\geq1\) and both
  \(v\in\{r,r+1\}\).
\end{enumerate}
The set \(\mathcal E\) is given explicitly in
\eqref{eq:seed-exceptions}.
\end{theorem}

\begin{proof}
Treat \(a\) as a formal variable in \(A_a(X)\), and put
\[
  c(a)=(1-a)^6(3a^2-3a+1)^2.
\]
We have the polynomial identities
\[
  \det W(A_a(G_*))=-8c(a),\qquad
  \det W(A_a(\overline{G_*}))=8c(a).
\]
The first is a direct calculation. For the second, with \(J\)
the all-ones matrix of order \(6\),
\[
  A_a(\overline{G_*})
  =-A_a(G_*)+(6a-1)I_6+(1-a)J.
\]
Since \(Jz=(\boldsymbol1^{\mathsf T}z)\boldsymbol1\), induction
shows that the two walk matrices are related by a triangular column
transformation with diagonal entries \((-1)^j\), \(j=0,\ldots,5\),
and determinant \((-1)^{0+1+\cdots+5}=-1\).
Because \(3a^2-3a+1=3(a-1/2)^2+1/4\), both graphs are
\(A_{\alpha}\)-controllable for every \(0\leq\alpha<1\).

It remains to check the three spectral exclusions. Define
\begin{align*}
  d_X(a)&=\det A_a(X),\\
  u_X(a)&=\det\bigl(aA_a(X)+(2a-1)I_6\bigr),\\
  v_X(a)&=\det\bigl(aA_a(X)+(1-2a)I_6\bigr).
\end{align*}
For \(a>0\), their nonvanishing excludes, respectively,
\(0\), \((1-2a)/a\), and \((2a-1)/a\) from the base spectrum.
At \(a=0\), \(u_X(0)=v_X(0)=1\). Direct calculation gives
\begin{align*}
  d_{G_*}(a)&=(4a-1)(a^2+2a-1)(9a^2-2a-1),\\
  u_{G_*}(a)&=(a^2+2a-1)\nu(a),\qquad v_{G_*}(a)=c(a),
\end{align*}
where
\[
\begin{aligned}
  \nu(a)={}&72a^9+565a^8+86a^7-629a^6+34a^5\\
         &+228a^4-54a^3-23a^2+10a-1.
\end{aligned}
\]
Set
\[
\begin{aligned}
  P(a)&=(4a-1)(a^2+2a-1)(9a^2-2a-1)\nu(a),\\
  \widehat P(a)&=d_{\overline{G_*}}(a)u_{\overline{G_*}}(a)
          v_{\overline{G_*}}(a).
\end{aligned}
\]
Since \(c(\alpha)>0\) and \(P\) includes the common factor
\(a^2+2a-1\) of \(d_{G_*}\) and \(u_{G_*}\) once, the three determinants
for \(G_*\) are nonzero exactly when \(P(\alpha)\neq0\). Define
\begin{equation}\label{eq:seed-exceptions}
  \mathcal E=\{\alpha\in[0,1):P(\alpha)=0\}.
\end{equation}
This set is finite, since \(P(0)=1\).
By its factorization, \(d_{G_*}(\alpha)\neq0\) on \([0,1/5)\).
For \(0<\alpha<1/5\), put \(x=1/\alpha>5\).
Scaling the walk and matrix determinants shows that
\(F_{G_*}(x)g_{G_*}(x)\neq0\).
Since \(\Delta(G_*)=3\), Theorem~\ref{thm:eventual-reciprocal} applies.
Applying Corollary~\ref{cor:even-center-controllable} at \(r=3\)
then gives \(P(\alpha)\neq0\).
Thus \(\mathcal E\cap[0,1/5)=\varnothing\).

Outside \(\mathcal E\), \(G_*\) satisfies
Corollary~\ref{cor:even-center-controllable} for every \(r\).
To handle the exceptional set, form the Sylvester matrix from
\(27\) shifted coefficient rows of \(P\), followed by \(14\) of \(\widehat P\).
Its \(41\times41\) determinant, computed modulo \(11\), gives
\[
  \Res_a(P,\widehat P)\equiv5\pmod{11}.
\]
This nonzero resultant shows that \(P\) and \(\widehat P\) have no common zero,
so \(\widehat P(\alpha)\neq0\) on \(\mathcal E\).
The complement therefore satisfies the same criterion there
for every \(r\) and both central roots.
\end{proof}

For each parameter, these products form an infinite family of orders
\(12r\), with the path rooted at an internal vertex for \(r\geq2\).
Two is the minimum number of fixed base graphs for uniform coverage.
Indeed, if \(G\) has an edge, then
\(A(G)\) is nonzero with trace zero, so \(\lambda_{\min}(A(G))<0\);
since \(A_{1/2}(G)=Q(G)/2\) is positive semidefinite, continuity
gives singularity at some \(\alpha\in(0,1/2]\), violating
Corollary~\ref{cor:even-center-controllable}.
For an edgeless base, the product at \(r=1\) is a disjoint union of
edges whose walk columns are all \(\boldsymbol1\), so it is uncontrollable.

\section{An algebraic criterion for base graphs}
\label{sec:avoidance}

For general rooted graphs, algebraic conjugates can simplify the
spectral test in Corollary~\ref{cor:rooted-general}.
Fix \(1/2\leq\alpha<1\) and set \(K_{\alpha}=\Q(\alpha)\).
For an element algebraic over \(K_{\alpha}\), its minimal polynomial
is the monic irreducible polynomial in \(K_{\alpha}[t]\) vanishing at it;
the complex zeros of this polynomial are its algebraic conjugates
over \(K_{\alpha}\).
Both \(\phi(A_{\alpha}(G);t)\) and \(h_{\alpha,H,v}(t)\) belong to
\(K_{\alpha}[t]\), and
\[
A_{\alpha}(G)=(1-\alpha)Q(G)+(2\alpha-1)D(G)
\]
is positive semidefinite. Hence an irreducible factor with a negative
or nonreal zero can be omitted entirely from the spectral test.
To apply the following proposition, factor \(h_{\alpha,H,v}\) over
\(K_{\alpha}\), collect all such factors with their multiplicities into
\(b\), and let \(a=h_{\alpha,H,v}/b\). The spectral test then involves
only \(a\).

\begin{proposition}
\label{prop:alpha-algebraic-sufficient}
Let \(1/2\leq\alpha<1\), and let \(H^{(v)}\) be a rooted graph such that
\(v\) is an \(A_{\alpha}\)-Wronskian vertex. Suppose that, in
\(K_{\alpha}[t]\),
\[
  h_{\alpha,H,v}(t)=a(t)b(t),
\]
where \(a,b\in K_{\alpha}[t]\) are nonzero and every nonconstant
irreducible factor of \(b\) has a negative or nonreal zero. Then, for every
graph \(G\), the rooted product \(G\circ H^{(v)}\) is
\(A_{\alpha}\)-controllable if and only if \(G\) is
\(A_{\alpha}\)-controllable and \(a(\lambda)\neq0\) for every
\(\lambda\in\Spec(A_{\alpha}(G))\).
\end{proposition}

\begin{proof}
First suppose that both graphs have at least two vertices.
By Corollary~\ref{cor:rooted-general}, it suffices to show that
\(b\) does not vanish on the spectrum of \(A_{\alpha}(G)\).
Suppose that \(b(\lambda)=0\) at an eigenvalue \(\lambda\), and
choose an irreducible factor \(f\) of \(b\) with \(f(\lambda)=0\).
Up to a nonzero scalar, \(f\) is the minimal polynomial of
\(\lambda\) over \(K_{\alpha}\), so it divides
\(\phi(A_{\alpha}(G);t)\). Every zero of \(f\) is therefore
real and nonnegative, contrary to the hypothesis on \(b\).
Hence \(b\) has no zero in the base spectrum, and the criterion
reduces to the stated conditions on \(G\) and \(a\).

For the one-vertex cases, let \(K_1\) denote the graph with one vertex
and no edges. If \(H=K_1\), then \(h_{\alpha,H,v}=1\) and
\(G\circ H^{(v)}=G\), so the assertion is immediate.
If \(G=K_1\), then \(G\circ H^{(v)}=H\), the base matrix is
\([0]\), and \(G\) is \(A_{\alpha}\)-controllable.
The hypothesis on \(b\) gives \(b(0)\neq0\), since otherwise \(t\)
would be an irreducible factor with no negative or nonreal zero.
Thus \(H\) is \(A_{\alpha}\)-controllable exactly when
\(h_{\alpha,H,v}(0)=a(0)b(0)\neq0\), equivalently, when \(a(0)\neq0\).
\end{proof}

For example, take \(\alpha=1/2\) and root \(P_7\) at vertex \(2\).
Since \(\gcd(3,7)=1\), this is a Wronskian root by
Corollary~\ref{cor:rational-classification}.
Successive multiplication by \(A_{1/2}(P_7)+tE_{22}\) gives the walk
columns, whose determinant is
\[
  h_{1/2,P_7,2}(t)=-2^{-15}t^3(t^3-12t-8).
\]
The cubic has roots in \((-4,-2)\), \((-2,0)\), and \((2,4)\).
Substituting \(t=2s\) gives \(8(s^3-3s-1)\); the rational-root test
gives irreducibility over \(\Q\).
The proposition removes the entire cubic, including its positive root.
Thus \(G\circ P_7^{(2)}\) is \(Q\)-controllable exactly when \(G\) is
\(Q\)-controllable and \(Q(G)\) is nonsingular.

For rational \(\alpha=p/q\geq1/2\), we can work with integer polynomials.
Take \(p/q\) in lowest terms with \(q>0\), and set
\[
  \widetilde A(X)=pD(X)+(q-p)A(X)=qA_{\alpha}(X),
  \qquad
  \widetilde h_{H,v}(z)
  =\det W\bigl(\widetilde A(H)+zE_{vv}\bigr).
\]
If \(|V(H)|=m\) and \(|V(G)|=n\), then
\[
  h_{\alpha,H,v}(t)
  =q^{-\binom m2}\widetilde h_{H,v}(qt),
  \qquad
  \phi(A_{\alpha}(G);t)
  =q^{-n}\phi(\widetilde A(G);qt).
\]
After the substitution \(z=qt\), these identities allow the
factor condition in Proposition~\ref{prop:alpha-algebraic-sufficient}
to be checked with polynomials in \(\mathbb Z[z]\),
up to nonzero rational scalars.

Finally, for \(0\leq\alpha<1/2\), the base matrix is real symmetric but
may have negative eigenvalues. The same proof applies if every
nonconstant irreducible factor of \(b\) has a nonreal zero.
A negative zero alone does not exclude such a factor.

\section*{Funding}

This work was partially supported by the National Natural Science Foundation
of China [grant number 12271182].
The funder had no role in the research design, analysis, manuscript
preparation, or decision to submit the paper for publication.

\section*{Declaration of competing interest}

The authors declare that they have no known competing financial interests or
personal relationships that could have appeared to influence the work
reported in this paper.

\section*{Declaration of generative AI and AI-assisted technologies
in the manuscript preparation process}

The authors used OpenAI Codex to assist with language editing, proof
review and revision, and the formalization of selected arguments in Lean.
The authors reviewed and revised the AI-assisted material and take full
responsibility for the manuscript.

\begingroup
\small
\setlength{\bibsep}{1pt}
\bibliographystyle{abbrvnat} 
\bibliography{references}
\endgroup
\end{document}